\documentclass[12pt]{article}

\usepackage{amsmath, amssymb, amsthm}
\usepackage{hyperref}
\usepackage{fullpage}

\title{Estimating $\pi$ with a Coin}

\author{Jim Propp}

\date{November 24, 2025; revised March 10, 2026}

\begin{document}

\maketitle

\begin{abstract}
\noindent
We describe a simple Monte Carlo method for estimating $\pi$
by tossing a coin. Although the underlying Catalan-number series
identities appear implicitly in the probability theory literature,
the interpretation of $\frac{\pi}{4}$ presented here seems to be new.
\end{abstract}

Buffon's needle gives a well-known way to estimate $\pi$:
when a needle of length $L$ is repeatedly dropped onto 
a plane ruled with parallel lines spaced $L$ units apart,
the proportion of drops that result in the needle landing 
on one of the lines, under the hypothesis that the drops 
are random and independent, is $\frac{2}{\pi}$.

Here we describe a different way you can use randomness 
to estimate $\pi$. Toss a coin until the first time the 
observed proportion of heads exceeds $\frac12$ (that
is, until the cumulative number of heads exceeds the 
cumulative number of tails) and record that fraction;
for instance, if the tosses proceeded as tails, heads, 
tails, heads, heads you would record the fraction 
$\frac35$. Now perform a second run of coin-tosses,
starting afresh and again recording a fraction that 
exceeds $\frac12$ (starting afresh means that you ignore
all earlier tosses). Repeat. As you perform more and 
more such trials, the average of the fractions will 
approach $\frac{\pi}{4}$.

We will prove this result using methods from the theory 
of random walks and stopping times. Let $\tau$ be the 
(random) number of coin tosses required until heads 
outnumber tails for the first time. Letting $H_\tau$ be 
the number of heads seen up to time $\tau$, we will prove 
that the expected value of the proportion ${H_\tau}/{\tau}$ 
is $\frac{\pi}{4}$.


The related fact that the expected value of ${1}/{\tau}$
is $\frac{\pi}{2} - 1$ is not new: it appears implicitly 
in the Catalan-number and arcsine-series machinery found in 
classical sources such as \cite{Feller}.
What appears to be new is the coin-tossing interpretation 
of $\frac{\pi}{4}$ as the expected stopping-time value of
the proportion ${H_\tau}/{\tau}$.


Let $H_n$ and $T_n$ be the numbers of heads and tails in the 
first $n$ independent tosses of a fair coin, and let
\[
S_n = H_n - T_n,
\]
so that $(S_n)$ is a simple symmetric random walk on $\mathbb Z$ 
with $S_0=0$. 
It is well known that this walk is recurrent, returning to 0
infinitely often with probability 1. 
Since the walk is recurrent and has unit step size, 
it will almost surely visit 1, so we define the stopping time
\[
\tau = \min\{n\ge 1 : S_n > 0\},
\]
the first time the walk hits $+1$.

Because $S_n$ changes by $\pm1$ at each step, $\tau$ can 
only take odd values. For $\tau=2k+1$ to occur, the walk 
must stay $\le 0$ for the first $2k$ steps, be at $0$ at time 
$2k$, and then step from $0$ to $+1$. The number of such paths 
is the Catalan number $C_{k} = {2k \choose k}/(k+1)$; since each 
path of length $2k+1$ has probability $2^{-(2k+1)}$, we have
\[
\mathbb P(\tau = 2k+1) = \frac{1}{2 \cdot 4^k} \frac{1}{k+1} {2k \choose k},
\qquad k=0,1,2,\dots.
\]
When the stopping time $\tau$ is $2k+1$, we have $H_\tau = k+1$, $T_\tau = k$,
and $\frac{H_\tau}{\tau} = \frac{k+1}{2k+1}$.
Thus
\begin{eqnarray*}
\mathbb{E}(H_\tau/\tau) 
& = & \sum_{k \geq 0} \mathbb{P} (\tau = 2k+1) \ \frac{k+1}{2k+1} \\
& = & \sum_{k \geq 0} \frac{1}{2 \cdot 4^k} \frac{1}{k+1} {2k \choose k} \frac{k+1}{2k+1} \\
& = & \frac12 \sum_{k \geq 0} \frac{1}{4^k} \frac{1}{2k+1} {2k \choose k}
\end{eqnarray*}
But the power series for arcsine is
\[\arcsin(x) = 
\sum_{n = 0}^\infty \frac{1}{2^{2n}} {2n \choose n} \frac{x^{2n+1}}{2n+1},\]
so we obtain
\[\mathbb{E}(H_\tau/\tau) = \frac12 \arcsin(1) = \frac{\pi}{4}\]
as claimed.

\bigskip

If you are in a classroom or a math club meeting, you
can use parallelization, letting students toss their 
own coins and record their own fractions. Collating 
all those fractions will result in a larger sample size 
and a more accurate estimate, but you should not expect
much accuracy, since the error will decrease like
$1/N^{1/4}$ where $N$ is the number of coin-tosses
(see~\cite{CM} for relevant estimates).

Upon learning of our result, Matt Parker used 10,000 coin flips
that he had recorded years earlier to estimate $\pi$, and obtained 
3.2266. An error on the order of $1/(10{,}000)^{1/4} = 1/10$ is 
indeed roughly what you should expect. If you want to obtain 
something nearer to 3.14, a trillion coin flips might be 
required---a task that, at a rate of one coin flip per second, 
would take over 30,000 years.

Stefan Gerhold brought the 2025 preprint~\cite{GH} to
our attention. That paper is concerned with how, in a
highly idealized setting, families' sex-ratios might be
affected by different sex-based stopping-rules families might 
employ. In the paragraph that contains their equation (3.7) 
Gerhold and Hubalek prove a formula that is equivalent to ours.

One of the referees asked: ``Why not also a short comment at least 
on the effect of a surplus of 2, for instance?'' Excellent idea!
It turns out that if you toss a coin until the number of heads
first exceeds the number of tails by $+2$ rather than $+1$, then the
expected proportion of heads is $\ln 2$ rather than $\pi/4$. Replacing
$+2$ by a generic positive integer $a$ appears to give rise to expressions 
involving $\pi$ when $a$ is odd and $\ln 2$ when $a$ is even,
but we have not checked ChatGPT's derivations and therefore regard
its claims and formulas as conjectural. We leave further investigation
to others.

A very general framework for approximating $\pi$ using a coin
was explored by Flajolet, Pelletier, and Soria~\cite{FPS}.
We expect that some of their coin-based algorithms for
estimating $\pi$, such as {\tt Rama()}, are much more 
efficient than ours.

It should be noted that the fraction $H_\tau/\tau$
equals 1 half the time and exceeds 1/2 the rest of the time, 
implying that $\mathbb{E}[H_\tau/\tau]$ is greater than 3/4. 
It is even easier to see that $\mathbb{E}[H_\tau/\tau]$ must 
be less than 1. Hence the proof of the formula 
$\mathbb{E}[H_\tau/\tau] = \pi/4$ can be regarded as 
a probabilistic proof of the well-known fact
that $\pi$ lies strictly between 3 and 4.

\bigskip

Acknowledgments:
Thanks to Eric Severson for suggesting a simplification in our proof
and thanks to Ven Popov for catching bibliographic errors.
Thanks also the referees for their careful reading and helpful 
suggestions.

\end{document}